\documentclass[11pt]{article}
\usepackage{mathrsfs}
\usepackage{latexsym,amsmath,amssymb,amsfonts}
\usepackage{color,colordvi,amscd,indentfirst,graphics}
\allowdisplaybreaks

\newtheorem{theo}{Theorem}[section]

\newtheorem{lem}[theo]{Lemma}
\newtheorem{coro}[theo]{Corollary}

\newtheorem{claim}[theo]{Claim}

\def\qed{\hfill \rule{4pt}{7pt}}
\def\pf{\noindent {\it Proof.} }

\begin{document}

\title{Spanning trees with at most $4$ leaves in $K_{1,5}$-free graphs}
\author{Yuan Chen\footnote{E-mail address: chenyuanmath@hotmail.com.}\\
School of Mathematics and Computer Science\\
Wuhan Textile University\\
1 Fangzhi Road, Wuhan 430073, P.R. China\\
\medskip\\
Pham Hoang Ha\footnote{E-mail address: ha.ph@hnue.edu.vn (Corresponding author).}\\
Department of Mathematics\\
Hanoi National University of Education\\
136 XuanThuy Street, Hanoi, Vietnam\\
\medskip\\
Dang Dinh Hanh\footnote{E-mail address: ddhanhdhsphn@gmail.com.}\\
Department of Mathematics\\
Hanoi Architectural University\\
Km10 NguyenTrai Street, Hanoi, Vietnam\\}
\date{}

\maketitle{}

\bigskip

\begin{abstract}
In 2009, Kyaw proved that every $n$-vertex connected $K_{1,4}$-free
graph $G$ with $\sigma_4(G)\geq n-1$ contains a spanning tree with
at most $3$ leaves. In this paper, we prove an analogue of Kyaw's
result for connected $K_{1,5}$-free graphs. We show that every
$n$-vertex connected $K_{1,5}$-free graph $G$ with $\sigma_5(G)\geq
n-1$ contains a spanning tree with at most $4$ leaves. Moreover, the
degree sum condition ``$\sigma_5(G)\geq n-1$" is best possible.
\end{abstract}



\newpage

\section{Introduction}

In this paper, we only consider finite simple graphs. Let $G$ be a
graph with vertex set $V(G)$ and edge set $E(G)$. For any vertex
$v\in V(G)$, we use $N_G(v)$ and $d_G(v)$ (or $N(v)$ and $d(v)$ if
there is no ambiguity) to denote the set of neighbors of $v$ and the
degree of $v$ in $G$, respectively. For any $X\subseteq V(G)$, we
denote by $|X|$ the cardinality of $X$. We define
$N(X)=\bigcup\limits_{x\in X}N(x)$ and $d(X)=\sum\limits_{x\in
X}d(x)$. For $k\geq 1$, we let $N_k(X)=\{x\in V(G)\;|\;|N(x)\cap
X|=k\}$ and $N_{\geq k}(X)=\{x\in V(G)\;|\;|N(x)\cap X|\geq k\}$. We
use $G-X$ to denote the graph obtained from $G$ by deleting the
vertices in $X$ together with their incident edges. The subgraph of
$G$ induced by $X$ is denoted by $G[X]$. We define $G-uv$ to be the
graph obtained from $G$ by deleting the edge $uv\in E(G)$, and
$G+uv$ to be the graph obtained from $G$ by adding an edge $uv$
between two non-adjacent vertices $u$ and $v$ of $G$. We write $A:=
B$ to rename $B$ as $A$.

A subset $X\subseteq V(G)$ is called an \emph{independent set} of
$G$ if no two vertices of $X$ are adjacent in $G$.  The maximum size
of an independent set in $G$ is denoted by $\alpha(G)$. For $k\geq
1$, we define $\sigma_k(G)=\min\{\sum\limits_{i=1}^k
d(v_i)\;|\;\{v_1,\ldots,v_k\}$ is an independent set in $G$\}. For
$r\geq 1$, a graph is said to be \emph{$K_{1,r}$-free} if it does
not contain $K_{1,r}$ as an induced subgraph. A $K_{1,3}$-free graph
is also called a \emph{claw-free} graph. We use $K_n$ to denote the
complete graph on $n$ vertices.

Let $T$ be a tree. A vertex of degree one is a \emph{leaf} of $T$
and a vertex of degree at least three is a \emph{branch vertex} of
$T$. For two distinct vertices $u,v$ of $T$, we denote by $P_T[u,v]$
the unique path in $T$ connecting $u$ and $v$ and denote by
$d_T[u,v]$ the distance between $u$ and $v$ in $T$. We define the
\emph{orientation} of $P_T[u,v]$ is from $u$ to $v$. We refer
to~\cite{Di05} for terminology and notation not defined here.

There are several well-known conditions (such as the independence
number conditions and the degree sum conditions) ensuring that a
graph $G$ contains a spanning tree with a bounded number of leaves
or branch vertices (see the survey paper~\cite{OY11} and the
references cited therein for details). Win~\cite{Wi79} obtained a
sufficient condition related to the independence number for
$k$-connected graphs, which confirms a conjecture of Las
Vergnas~\cite{LV71}. Broersma and Tuinstra~\cite{BT98} gave a degree
sum condition for a connected graph to contain a spanning tree with
at most $m$ leaves.

\begin{theo}\label{theo1.1}{\rm (Win~\cite{Wi79})}
Let $G$ be a $k$-connected graph and let $m\geq 2$.  If
$\alpha(G)\leq k+m-1$, then $G$ has a spanning tree with at most $m$
leaves.
\end{theo}

\begin{theo}\label{theo1.2}{\rm (Broerma and Tuinstra~\cite{BT98})}
Let $G$ be a connected graph with $n$ vertices and let $m\geq 2$. If
$\sigma_2(G)\geq n-m+1$, then $G$ has a spanning tree with at most
$m$ leaves.
\end{theo}

Kano et al.~\cite{KKMOSY12} presented a degree sum condition for a
connected claw-free graph to have a spanning tree with at most $m$
leaves, which generalizes a result of Matthews and
Sumner~\cite{MS84} and a result of Gargano et al.~\cite{GHHSV04}.
Matsuda, Ozeki and Yamashita~\cite{MOY14} and Chen, Li and
Xu~\cite{CLX17} considered the sufficient conditions for a connected
claw-free graph to have a spanning tree with few branch vertices or
few leaves, respectively.

\begin{theo}\label{theo1.3}{\rm (Kano et al.~\cite{KKMOSY12})}
Let $G$ be a connected claw-free graph with $n$ vertices and let
$m\geq 2$. If $\sigma_{m+1}(G)\geq n-m$, then $G$ has a spanning
tree with at most $m$ leaves.
\end{theo}

\begin{theo}\label{theo1.4}{\rm (Matsuda, Ozeki and Yamashita~\cite{MOY14})}
Let $G$ be a connected claw-free graph with $n$ vertices. If
$\sigma_5(G)\geq n-2$, then $G$ contains a spanning tree with at
most one branch vertex.
\end{theo}

\begin{theo}\label{theo1.5}{\rm (Chen, Li and Xu~\cite{CLX17})}
Let $G$ be a $k$-connected claw-free graph with $n$ vertices. If
$\sigma_{k+3}(G)\geq n-k$, then $G$ contains a spanning tree with at
most $3$ leaves.
\end{theo}

For connected $K_{1,4}$-free graphs, Kyaw~\cite{Ky09,Ky11} obtained
the following two sharp results.

\begin{theo}\label{theo1.6}{\rm (Kyaw~\cite{Ky09})}
Let $G$ be a connected $K_{1,4}$-free graph with $n$ vertices. If
$\sigma_4(G)\geq n-1$, then $G$ contains a spanning tree with at
most $3$ leaves.
\end{theo}

\begin{theo}\label{theo1.7}{\rm (Kyaw~\cite{Ky11})}
Let $G$ be a connected $K_{1,4}$-free graph with $n$ vertices.
\begin{itemize}
\item [$($i$)$]
If $\sigma_3(G)\geq n$, then $G$ has a hamiltonian path.
\item [$($ii$)$]
If $\sigma_{m+1}(G)\geq n-\frac{m}{2}$ for some integer $m\geq 3$,
then $G$ has a spanning tree with at most $m$ leaves.
\end{itemize}
\end{theo}

Chen, Chen and Hu~\cite{CCH14} considered the degree sum condition
for a $k$-connected $K_{1,4}$-free graph to contain a spanning tree
with at most $3$ leaves.

\begin{theo}\label{theo1.8}{\rm (Chen, Chen and Hu~\cite{CCH14})}
Let $G$ be a $k$-connected $K_{1,4}$-free graph with $n$ vertices
and let $k\geq2$. If $\sigma_{k+3}(G)\geq n+2k-2$, then $G$ has a
spanning tree with at most $3$ leaves.
\end{theo}

In this paper, we further consider connected $K_{1,5}$-free graphs.
We give a sufficient condition for a connected $K_{1,5}$-free graph
to have a spanning tree with few leaves.

\begin{theo}\label{theo1.9}
Let $G$ be a connected $K_{1,5}$-free graph with $n$ vertices. If
$\sigma_5(G)\geq n-1$, then $G$ contains a spanning tree with at
most $4$ leaves.
\end{theo}

It is easy to see that if a tree has at most $k$ leaves ($k\geq 2$),
then it has at most $k-2$ branch vertices. Therefore, we immediately
obtain the following corollary from Theorem~\ref{theo1.9}.

\begin{coro}\label{coro1.10}
Let $G$ be a connected $K_{1,5}$-free graph with $n$ vertices. If
$\sigma_5(G)\geq n-1$, then $G$ contains a spanning tree with at
most $2$ branch vertices.
\end{coro}

We end this section by constructing an example to show that the
degree sum condition ``$\sigma_5(G)\geq n-1$" in
Theorem~\ref{theo1.9} is sharp. For $m\geq 1$, let
$D_1,D_2,D_3,D_4,D_5$ be vertex-disjoint copies of $K_m$ and let
$xy$ be an edge such that neither $x$ nor $y$ is contained in
$\bigcup\limits_{i=1}^5V(D_i)$. Join $x$ to all the vertices in
$D_1\cup D_2\cup D_3$ and join $y$ to all the vertices in $D_4\cup
D_5$. The resulting graph is denoted by $G$. Then it is easy to
check that $G$ is a connected $K_{1,5}$-free graph with $n=5m+2$
vertices and $\sigma_5(G)=5m=n-2$. However, every spanning tree of
$G$ contains at least $5$ leaves.

\section{Proof of the main result}

In this section, we extend the idea of Kyaw in~\cite{Ky09} to prove
Theorem~\ref{theo1.9}. For this purpose, we need the following
lemma.

\begin{lem}\label{lem2.1}
Let $G$ be a connected graph such that $G$ does not have a spanning
tree with at most $4$ leaves, and let $T$ be a maximal tree of $G$
with $5$ leaves. Then there does not exist a tree $T'$ in $G$ such
that $T'$ has at most $4$ leaves and $V(T')=V(T)$.
\end{lem}

\pf Suppose for a contradiction that there exists a tree $T'$ in $G$
with at most $4$ leaves and $V(T')=V(T)$. Since $G$ has no spanning
tree with at most $4$ leaves, we see that $V(G)-V(T')\neq\emptyset$.
Hence there must exist two vertices $v$ and $w$ in $G$ such that
$v\in V(T')$ and $w\in N(v)\cap (V(G)-V(T'))$. Let $T_1$ be the tree
obtained from $T'$ by adding the vertex $w$ and the edge $vw$.

If $T_1$ has $5$ leaves, then $T_1$ contradicts the maximality of
$T$ (since $|V(T_1)|=|V(T)|+1>|V(T)|$). So we may assume that $T_1$
has at most $4$ leaves. By repeating this process, we can
recursively construct a set of trees $\{T_i\;|\;i\geq 1\}$ in $G$
such that $T_i$ has at most $4$ leaves and $|V(T_{i+1})|=|V(T_i)|+1$
for each $i\geq 1$. Since $G$ has no spanning tree with at most $4$
leaves and $|V(G)|$ is finite, the process must terminate after a
finite number of steps, i.e., there exists some $k\geq 1$ such that
$T_{k+1}$ is a tree in $T$ with $5$ leaves. But this contradicts the
maximality of $T$. So the lemma holds.    \qed

\bigskip

\noindent{\bf Proof of Theorem~\ref{theo1.9}.} We prove the theorem
by contradiction. Suppose to the contrary that $G$ contains no
spanning tree with at most $4$ leaves. Then every spanning tree of
$G$ contains at least $5$ leaves. We choose a maximal tree $T$ of
$G$ with exactly $5$ leaves. Let $U=\{u_1,u_2,u_3,u_4,u_5\}$ be the
set of leaves of $T$. By the maximality of $T$, we have
$N(U)\subseteq V(T)$.

We consider three cases according to the number of branch vertices
in $T$. (Note that $T$ contains at most three branch vertices.)

\medskip

\textit{Case 1}. $T$ contains two branch vertices.

\medskip

Let $s$ and $t$ be the two branch vertices in $T$ such that
$d_T(s)=4$ and $d_T(t)=3$. For each $1\leq i\leq 5$, let $B_i$ be
the vertex set of the connected component of $T-\{s,t\}$ containing
$u_i$ and let $v_i$ be the unique vertex in $B_i\cap N_T(\{s,t\})$.
Without loss of generality, we may assume that
$\{v_1,v_2,v_3\}\subseteq N_T(s)$ and $\{v_4,v_5\}\subseteq N_T(t)$.
For each $1\leq i\leq 5$ and $x\in B_i$, we use $x^-$ and $x^+$ to
denote the predecessor and the successor of $x$ on $P_T[s,u_i]$ or
$P_T[t,u_i]$, respectively (if such a vertex exists). Let $s^+$ and
$t^-$ be the successor of $s$ and the predecessor of $t$ on
$P_{T}[s,t]$, respectively. Define $P:=V(P_T[s,t])-\{s,t\}$.

For this case, we further choose $T$ such that
\begin{itemize}
\item [$($C1$)$]
$d_T[s,t]$ is as small as possible, and
\item [$($C2$)$]
subject to (C1), $\sum\limits_{i=1}^3 |B_i|$ is as large as
possible.
\end{itemize}

\begin{claim}\label{claim2.2}
For all $1\leq i,j\leq 5$ and $i\neq j$, if $x\in N(u_j)\cap B_i$,
then $x\notin \{u_i,v_i\}$ and $x^-\notin N(U-\{u_j\})$.
\end{claim}

\pf Suppose $x\in \{u_i,v_i\}$. Then $T':=T-v_iv_i^-+xu_j$ is a tree
in $G$ with $4$ leaves and $V(T')=V(T)$, which contradicts
Lemma~\ref{lem2.1}. So we have $x\notin \{u_i,v_i\}$.

Next, assume $x^-\in N(U-\{u_j\})$. Then there exists some
$k\in\{1,2,3,4,5\}-\{j\}$ such that $x^-u_k\in E(G)$. Now,
$T':=T-\{v_iv_i^-,xx^-\}+\{xu_j,x^-u_k\}$ is a tree in $G$ with $4$
leaves and $V(T')=V(T)$, also contradicting Lemma~\ref{lem2.1}. This
proves Claim~\ref{claim2.2}.     \qed

\bigskip

By Claim~\ref{claim2.2}, we know that $U$ is an independent set in
$G$. Since $G$ is $K_{1,5}$-free, we have $N_5(U)=\emptyset$.

\begin{claim}\label{claim2.3}
$N(u_i)\cap P=\emptyset$ for each $4\leq i\leq 5$.
\end{claim}

\pf Suppose the assertion of the claim is false. Then there exists
some vertex $x\in P$ such that $xu_i\in E(G)$ for some
$i\in\{4,5\}$. Let $T':=T-tv_i+xu_i$, then $T'$ is a tree in $G$
with $5$ leaves such that $V(T')=V(T)$, $T'$ has two branch vertices
$s$ and $x$, $d_{T'}(s)=4$, $d_{T'}(x)=3$ and
$d_{T'}[s,x]<d_T[s,t]$. But this contradicts the condition (C1). So
the claim holds.   \qed

\begin{claim}\label{claim2.4}
If $P\neq \emptyset$, then $\sum\limits_{i=1}^3 |N(u_i)\cap
\{x\}|\leq 1$ for each $x\in P$.
\end{claim}

\pf Suppose to the contrary that there exists some vertex $x\in P$
such that $\sum\limits_{i=1}^3 |N(u_i)\cap \{x\}|\geq 2$. Then there
exist two distinct $j,k\in\{1,2,3\}$ such that $xu_j,xu_k\in E(G)$.
Let $T':=T-\{sv_j,sv_k\}+\{xu_j,xu_k\}$, then $T'$ is a tree in $G$
with $5$ leaves such that $V(T')=V(T)$, $T'$ has two branch vertices
$x$ and $t$, $d_{T'}(x)=4$, $d_{T'}(t)=3$ and
$d_{T'}[x,t]<d_T[s,t]$, contradicting the condition (C1). This
completes the proof of Claim~\ref{claim2.4}.     \qed

\begin{claim}\label{claim2.5}
If $P\neq \emptyset$, then $N(U)\cap \{s^+\}=\emptyset$.
\end{claim}

\pf Suppose this is false. Then by Claim~\ref{claim2.3}, there
exists some $i\in\{1,2,3\}$ such that $s^+u_i\in E(G)$. Now,
$T':=T-ss^++s^+u_i$ is a tree in $G$ with $4$ leaves and
$V(T')=V(T)$, which contradicts Lemma~\ref{lem2.1}. So the assertion
of the claim holds.     \qed

\begin{claim}\label{claim2.6}
$N(u_i)\cap \{s\}=\emptyset$ for each $4\leq i\leq 5$.
\end{claim}

\pf Suppose $su_i\in E(G)$ for some $i\in\{4,5\}$. If $P=\emptyset$,
then we have $st\in E(T)$ and $T':=T-st+su_i$ is a tree in $G$ with
$4$ leaves and $V(T')=V(T)$, contradicting Lemma~\ref{lem2.1}. So we
may assume that $P\neq\emptyset$ and hence $s^+\neq t$. By applying
Claims~\ref{claim2.2} and~\ref{claim2.3}, we deduce that
$N(u_i)\cap\{s^+,v_1,v_2,v_3\}=\emptyset$.

Suppose that $s^+v_j\in E(G)$ for some $j\in\{1,2,3\}$. Then
$T':=T-\{ss^+,sv_j\}+\{su_i,s^+v_j\}$ is a tree in $G$ with $4$
leaves and $V(T')=V(T)$, which contradicts Lemma~\ref{lem2.1}. So we
conclude that $N(s^+)\cap \{v_1,v_2,v_3\}=\emptyset$.

Now, assume there exits two distinct $j,k\in \{1,2,3\}$ such that
$v_jv_k\in E(G)$. Then by Claim~\ref{claim2.2}, we see that $u_k\neq
v_k$. Let $T':=T-\{sv_j,tv_i\}+\{su_i,v_jv_k\}$, then $T'$ is a tree
in $G$ with $5$ leaves such that $V(T')=V(T)$, $T'$ has two branch
vertices $s$ and $v_k$, $d_{T'}(s)=4$, $d_{T'}(v_k)=3$ and
$d_{T'}[s,v_k]<d_T[s,t]$, contradicting the condition (C1).
Therefore, $v_1,v_2$ and $v_3$ are pairwise non-adjacent in $G$.

But then, $\{s^+,u_i,v_1,v_2,v_3\}$ is an independent set and
$G[\{s,s^+,u_i,v_1,v_2,v_3\}]$ is an induced $K_{1,5}$ of $G$, again
a contradiction. This proves Claim~\ref{claim2.6}.      \qed

\begin{claim}\label{claim2.7}
If $\sum\limits_{i=1}^5 |N(u_i)\cap \{t\}|\geq 3$, then
$P\neq\emptyset$.
\end{claim}

\pf Suppose for a contradiction that $P=\emptyset$. Then we have
$st\in E(G)$. Since $\sum\limits_{i=1}^5 |N(u_i)\cap \{t\}|\geq 3$,
there exists some $j\in\{1,2,3\}$ such that $tu_j\in E(G)$. Let
$T':=T-st+tu_j$, then $T'$ is a tree in $G$ with $4$ leaves and
$V(T')=V(T)$, which contradicts Lemma~\ref{lem2.1}. So the claim
holds.    \qed

\begin{claim}\label{claim2.8}
$N_4(U)=\emptyset$.
\end{claim}

\pf Suppose to the contrary that there exists some vertex $x\in
N_4(U)$. Then by Claims~\ref{claim2.3} and~\ref{claim2.6}, we have
$x\in B_1\cup B_2\cup B_3\cup B_4\cup B_5\cup\{t\}$.

First, suppose $x\in B_i$ for some $1\leq i\leq 5$. By
Claim~\ref{claim2.2}, we know that $x^-\notin N(U)$. Then $(N(x)\cap
U)\cup\{x^-\}$ is an independent set and $G[(N(x)\cap
U)\cup\{x,x^-\}]$ is an induced $K_{1,5}$ of $G$, contradicting the
assumption that $G$ is $K_{1,5}$-free.

So we may assume that $x=t$. Then by Claim~\ref{claim2.7}, we
conclude that $P\neq \emptyset$ and hence $t^-\neq s$. It follows
from Claim~\ref{claim2.3} that $N(u_i)\cap\{t^-\}=\emptyset$ for
each $4\leq i\leq 5$. Suppose that $t^-u_j\in E(G)$ for some
$j\in\{1,2,3\}$. Since $t\in N_4(U)$, there exists some
$k\in\{1,2,3\}-\{j\}$ such that $tu_k\in E(G)$. Let
$T':=T-\{sv_j,tt^-\}+\{tu_k,t^-u_j\}$, then $T'$ is a tree in $G$
with $4$ leaves and $V(T')=V(T)$, which contradicts
Lemma~\ref{lem2.1}. Therefore, we deduce that $N(U)\cap
\{t^-\}=\emptyset$. But then, $(N(t)\cap U)\cup\{t^-\}$ is an
independent set and $G[(N(t)\cap U)\cup\{t,t^-\}]$ is an induced
$K_{1,5}$ of $G$, again a contradiction. This completes the proof of
Claim~\ref{claim2.8}. \qed

\begin{claim}\label{claim2.9}
$(N_3(U)-N(u_i))\cap B_i=\emptyset$ for each $1\leq i\leq 5$.
\end{claim}

\pf Suppose this is false. Then there exists some vertex $x\in
(N_3(U)-N(u_i))\cap B_i$ for some $1\leq i\leq 5$. By applying
Claim~\ref{claim2.2}, we have $x\notin\{u_i,v_i\}$ and
$x^-,x^+\notin N(U-\{u_i\})$.

Suppose that $x^-x^+\in E(G)$. Since $x\in N_3(U)-N(u_i)$, there
must exist two distinct $j,k\in\{1,2,3,4,5\}-\{i\}$ such that
$xu_j,xu_k\in E(G)$. Then
$T':=T-\{v_jv_j^-,xx^-,xx^+\}+\{xu_j,xu_k,x^-x^+\}$ is a tree in $G$
with $4$ leaves and $V(T')=V(T)$, contradicting Lemma~\ref{lem2.1}.
Hence $x^-x^+\notin E(G)$.

Now, $(N(x)\cap U)\cup\{x^-,x^+\}$ is an independent set and
$G[(N(x)\cap U)\cup\{x,x^-,x^+\}]$ is an induced $K_{1,5}$ of $G$,
giving a contradiction. So the assertion of the claim holds.   \qed

\begin{claim}\label{claim2.10}
$N(u_j)\cap B_i=\emptyset$ for all $4\leq i\leq 5$ and $1\leq j\leq
3$. In particular, $N_3(U)\cap N(u_i)\cap B_i=\emptyset$ for each
$4\leq i\leq 5$.
\end{claim}

\pf Suppose the assertion of the claim is false. Then there exists
some vertex $x\in B_i$ such that $xu_j\in E(G)$ for some
$i\in\{4,5\}$ and $j\in\{1,2,3\}$. By Claim~\ref{claim2.2}, we have
$x\notin\{u_i,v_i\}$. Let $T':=T-xx^-+xu_j$, and let $B_k'$ be the
vertex set of the connected component of $T'-\{s,t\}$ containing
$u_k$ for each $1\leq k\leq 3$. It is easy to check that $T'$ is a
tree in $G$ with $5$ leaves such that $V(T')=V(T)$, $T'$ has two
branch vertices $s$ and $t$, $d_{T'}(s)=4$, $d_{T'}(t)=3$,
$d_{T'}[s,t]=d_T[s,t]$ and $\sum\limits_{k=1}^3
|B_k'|=\sum\limits_{k=1}^3 |B_k|+|V(P_T[x,u_i])|>\sum\limits_{k=1}^3
|B_k|$. But this contradicts the condition (C2). This proves
Claim~\ref{claim2.10}. \qed

\begin{claim}\label{claim2.11}
$|N_3(U)\cap N(u_i)\cap B_i|\leq 1$ for each $1\leq i\leq 3$.
\end{claim}

\pf Suppose for a contradiction that there exist two distinct
vertices $x,y\in N_3(U)\cap N(u_i)\cap B_i$ for some
$i\in\{1,2,3\}$. Without loss of generality, we may assume that
$x\in V(P_T[s,y])$. By Claim~\ref{claim2.2}, we have
$x,y\notin\{u_i,v_i\}$, $x^-\notin N(U)$ and $x^+\notin
N(U-\{u_i\})$. In particular, $x^+\neq y$. Since $x,y\in N_3(U)\cap
N(u_i)$, there exist two distinct $j,k\in\{1,2,3,4,5\}-\{i\}$ such
that $xu_j,yu_k\in E(G)$. We may assume that $x^-x^+,x^+u_i\notin
E(G)$; for otherwise,
$$T':=\left\{\begin{array}{ll}T-\{sv_i,xx^-,xx^+,yy^+\}+\{xu_i,xu_j,x^-x^+,yu_k\},
& \;\mbox{ if } x^-x^+\in E(G),\\
T-\{sv_i,xx^+,yy^-\}+\{xu_j,x^+u_i,yu_k\}, & \;\mbox{ if } x^+u_i\in
E(G),
\end{array}\right.$$
is a tree in $G$ with $4$ leaves and $V(T')=V(T)$, which contradicts
Lemma~\ref{lem2.1}. But then, $(N(x)\cap U)\cup\{x^-,x^+\}$ is an
independent set and $G[(N(x)\cap U)\cup\{x,x^-,x^+\}]$ is an induced
$K_{1,5}$ of $G$, again a contradiction. So the claim holds.
\qed

\begin{claim}\label{claim2.12}
For each $1\leq i\leq 3$, if $u_iv_i\in E(G)$, then $N_3(U)\cap
N(u_i)\cap B_i=\emptyset$.
\end{claim}

\pf Suppose to the contrary that $u_iv_i\in E(G)$ and there exists
some vertex $x\in N_3(U)\cap N(u_i)\cap B_i$ for some
$i\in\{1,2,3\}$. By Claim~\ref{claim2.2}, we have $x\neq v_i$. Since
$x\in N_3(U)\cap N(u_i)$, there exists some
$j\in\{1,2,3,4,5\}-\{i\}$ such that $xu_j\in E(G)$. Let
$T':=T-\{sv_i,xx^-\}+\{u_iv_i,xu_j\}$, then $T'$ is a tree in $G$
with $4$ leaves and $V(T')=V(T)$, contradicting Lemma~\ref{lem2.1}.
This completes the proof of Claim~\ref{claim2.12}.     \qed

\begin{claim}\label{claim2.13}
For each $1\leq i\leq 3$, if $su_i\in E(G)$, then $N_3(U)\cap
N(u_i)\cap B_i=\emptyset$.
\end{claim}

\pf For the sake of convenience, we may assume by symmetry that
$i=1$. Suppose the assertion of the claim is false. Then there
exists some vertex $x\in N_3(U)\cap N(u_1)\cap B_1$. By applying
Claims~\ref{claim2.2} and~\ref{claim2.12}, we know that
$x\notin\{u_1,v_1\}$ and $N(u_1)\cap \{v_1,v_2,v_3\}=\emptyset$.

Suppose $v_1v_j\in E(G)$ for some $j\in\{2,3\}$. Then
$T':=T-\{sv_1,sv_j\}+\{su_1,v_1v_j\}$ is a tree in $G$ with $4$
leaves and $V(T')=V(T)$, which contradicts Lemma~\ref{lem2.1}. So we
have $v_1v_2,v_1v_3\notin E(G)$.

Next, assume that $v_2v_3\in E(G)$. Then $u_2\neq v_2$ and $u_3\neq
v_3$ by Claim~\ref{claim2.2}. If there exists some $j\in\{2,3\}$
such that $xu_j\in E(G)$, then $T':=T-\{sv_2,sv_3\}+\{v_2v_3,xu_j\}$
is a tree in $G$ with $4$ leaves and $V(T')=V(T)$, contradicting
Lemma~\ref{lem2.1}. Hence $xu_2,xu_3\notin E(G)$. Since $x\in
N_3(U)\cap N(u_1)$, we conclude that $xu_4,xu_5\in E(G)$. Let
$T':=T-\{sv_2,tt^-,xx^-\}+\{su_1,v_2v_3,xu_4\}$. If $P=\emptyset$,
then $t^-=s$, and $T'$ is a tree in $G$ with $4$ leaves and
$V(T')=V(T)$, giving a contradiction. So we deduce that
$P\neq\emptyset$. But then, $T'$ is a tree in $G$ with $5$ leaves
such that $V(T')=V(T)$, $T'$ has two branch vertices $s$ and $v_3$,
$d_{T'}(s)=4$, $d_{T'}(v_3)=3$ and $d_{T'}[s,v_3]<d_T[s,t]$,
contradicting the condition (C1). Therefore, $v_1,v_2$ and $v_3$ are
pairwise non-adjacent in $G$.

We now consider the vertex $s^+$. We will show that
$N(s^+)\cap\{u_1,v_1,v_2,v_3\}=\emptyset$.

We first prove that $s^+u_1\notin E(G)$. Suppose this is false. Then
by Claim~\ref{claim2.5}, we see that $P=\emptyset$ and hence
$s^+=t$. Let $T':=T-st+tu_1$, then $T'$ is a tree in $G$ with $4$
leaves and $V(T')=V(T)$, which contradicts Lemma~\ref{lem2.1}.

We then have $s^+v_1\notin E(G)$; for otherwise,
$T':=T-\{ss^+,sv_1\}+\{su_1,s^+v_1\}$ is a tree in $G$ with $4$
leaves and $V(T')=V(T)$, also contradicting Lemma~\ref{lem2.1}.

Finally, we show that $s^+v_2,s^+v_3\notin E(G)$. Suppose not, and
let $s^+v_j\in E(G)$ for some $j\in\{2,3\}$. If there exists some
$k\in\{4,5\}$ such that $xu_k\in E(G)$, then
$T':=T-\{ss^+,sv_j\}+\{s^+v_j,xu_k\}$ is a tree in $G$ with $4$
leaves and $V(T')=V(T)$, which contradicts Lemma~\ref{lem2.1}.
Therefore, we have $xu_4,xu_5\notin E(G)$. Since $x\in N_3(U)\cap
N(u_1)$, we deduce that $xu_2,xu_3\in E(G)$. Let
$T':=T-\{ss^+,sv_j,xx^-,xx^+\}+\{su_1,s^+v_j,xu_2,xu_3\}$, then $T'$
is a tree in $G$ with $4$ leaves and $V(T')=V(T)$, again a
contradiction. Hence $N(s^+)\cap\{u_1,v_1,v_2,v_3\}=\emptyset$.

Now, $\{s^+,u_1,v_1,v_2,v_3\}$ is an independent set and
$G[\{s,s^+,u_1,v_1,v_2,v_3\}]$ is an induced $K_{1,5}$ of $G$,
giving a contradiction. So the assertion of the claim holds. \qed

\bigskip

By Claim~\ref{claim2.2}, $\{u_i\}$, $N(u_i)\cap B_i$,
$(N(U-\{u_i\}))^-\cap B_i$ and $(N_2(U)-N(u_i))\cap B_i$ are
pairwise disjoint subsets in $B_i$ for each $1\leq i\leq 5$, where
$(N(U-\{u_i\}))^-\cap B_i=\{x^-\;|\;x\in N(U-\{u_i\})\cap B_i\}$.
Recall that $N_5(U)=N_4(U)=(N_3(U)-N(u_i))\cap B_i=\emptyset$ (for
each $1\leq i\leq 5$) by Claims~\ref{claim2.8} and~\ref{claim2.9}.
Then for each $1\leq i\leq 3$, we conclude that
\begin{align*}
|B_i|&\geq 1+|N(u_i)\cap B_i|+|(N(U-\{u_i\}))^-\cap
B_i|+|(N_2(U)-N(u_i))\cap B_i|\\
&= 1+|N(u_i)\cap B_i|+|N(U-\{u_i\})\cap B_i|+|(N_2(U)-N(u_i))\cap
B_i|\\
&=1+\sum\limits_{j=1}^5|N(u_j)\cap B_i|-|N_3(U)\cap N(u_i)\cap
B_i|\\
&\geq\sum\limits_{j=1}^5|N(u_j)\cap B_i|+|N(u_i)\cap\{s\}|,
\tag{$1$}
\end{align*}
where the last inequality follows from Claims~\ref{claim2.11}
and~\ref{claim2.13}. Similarly, for each $4\leq i\leq 5$, we have
\begin{align*}
|B_i|&\geq 1+|N(u_i)\cap B_i|+|(N(U-\{u_i\}))^-\cap
B_i|+|(N_2(U)-N(u_i))\cap B_i|\\
&= 1+|N(u_i)\cap B_i|+|N(U-\{u_i\})\cap B_i|+|(N_2(U)-N(u_i))\cap
B_i|\\
&=1+\sum\limits_{j=1}^5|N(u_j)\cap B_i|-|N_3(U)\cap N(u_i)\cap
B_i|\\
&=1+\sum\limits_{j=1}^5|N(u_j)\cap B_i|+|N(u_i)\cap\{s\}|, \tag{$2$}
\end{align*}
where the last equality follows from Claims~\ref{claim2.6}
and~\ref{claim2.10}.

For each $1\leq i\leq 5$, we define $d_i=|N(u_i)\cap P|$. Then
$d_4=d_5=0$ by Claim~\ref{claim2.3}. By applying
Claim~\ref{claim2.4}, we know that $N(u_1)\cap P, N(u_2)\cap P$ and
$N(u_3)\cap P$ are pairwise disjoint. Therefore,
\begin{align*}
|P|\geq \sum\limits_{i=1}^5 d_i=\sum\limits_{i=1}^5|N(u_i)\cap P|.
\end{align*}
Moreover, if $P\neq\emptyset$, then by Claim~\ref{claim2.5}, we see
that $N(U)\cap\{s^+\}=\emptyset$ and hence
\begin{align*} |P|\geq 1+\sum\limits_{i=1}^5
d_i=1+\sum\limits_{i=1}^5|N(u_i)\cap P|.
\end{align*}
It follows from Claim~\ref{claim2.8} that
$\sum\limits_{i=1}^5|N(u_i)\cap \{t\}|\leq 3$. If
$\sum\limits_{i=1}^5|N(u_i)\cap \{t\}|\leq 2$, then
\begin{align*}
|V(P_T[s,t])|=2+|P|\geq
\sum\limits_{i=1}^5|N(u_i)\cap\{t\}|+\sum\limits_{i=1}^5|N(u_i)\cap
P|.
\end{align*}
If $\sum\limits_{i=1}^5|N(u_i)\cap \{t\}|=3$, then by
Claim~\ref{claim2.7}, we deduce that $P\neq \emptyset$. This implies
that
\begin{align*}
|V(P_T[s,t])|=2+|P|\geq 2+\left(1+\sum\limits_{i=1}^5|N(u_i)\cap
P|\right)=
\sum\limits_{i=1}^5|N(u_i)\cap\{t\}|+\sum\limits_{i=1}^5|N(u_i)\cap
P|.
\end{align*}
In both cases, we have
\begin{align*}
|V(P_T[s,t])|\geq
\sum\limits_{i=1}^5|N(u_i)\cap\{t\}|+\sum\limits_{i=1}^5|N(u_i)\cap
P|.   \tag{$3$}
\end{align*}

Note that $N(U)\subseteq V(T)$. By (1), (2) and (3), we conclude
that
\begin{align*}
|V(T)|&=\sum\limits_{i=1}^{3}|B_i|+\sum\limits_{i=4}^{5}|B_i|+|V(P_{T}[s,t])|\\
&\geq\sum\limits_{i=1}^{3}\left(\sum\limits_{j=1}^5|N(u_j)\cap
B_i|+|N(u_i)\cap\{s\}|\right)+
\sum\limits_{i=4}^{5}\left(1+\sum\limits_{j=1}^5|N(u_j)\cap
B_i|+|N(u_i)\cap\{s\}|\right)\\
&\ \ \ \
+\left(\sum\limits_{i=1}^5|N(u_i)\cap\{t\}|+\sum\limits_{i=1}^5|N(u_i)\cap
P|\right)\\
&=2+\sum\limits_{i=1}^{5}\sum\limits_{j=1}^{5}|N(u_j)\cap
B_i|+\sum\limits_{i=1}^{5}|N(u_i)\cap \{s,t\}|+\sum\limits_{i=1}^{5}|N(u_i)\cap P|\\
&=\sum\limits_{j=1}^{5}|N(u_j)\cap V(T)|+2\\
&=\sum\limits_{j=1}^{5}d(u_j)+2\\
&=d(U)+2.
\end{align*}
Since $U$ is an independent set in $G$, we have
\begin{align*}
n-1\leq \sigma_5(G)\leq d(U)\leq |V(T)|-2\leq n-2,
\end{align*}
a contradiction.

\medskip

\textit{Case 2}. $T$ contains only one branch vertex.

\medskip

Let $r$ be the unique branch vertex in $T$ with $d_T(r)=5$ and let
$N_T(r)=\{v_1,v_2,v_3,v_4,v_5\}$. Since $G$ is $K_{1,5}$-free, there
exist two distinct $i,j\in\{1,2,3,4,5\}$ such that $v_iv_j\in E(G)$.
Let $T':=T-rv_i+v_iv_j$. If $v_j$ is a leaf of $T$, then $T'$ is a
tree in $G$ with $4$ leaves and $V(T')=V(T)$, which contradicts
Lemma~\ref{lem2.1}. So we may assume that $v_j$ has degree two in
$T$. Then $T'$ is a tree in $G$ with $5$ leaves such that
$V(T')=V(T)$, $T'$ has two branch vertices $r$ and $v_j$,
$d_{T'}(r)=4$ and $d_{T'}(v_j)=3$. By the same argument as in the
proof of Case 1, we can also derive a contradiction.

\medskip

\textit{Case 3}. $T$ contains three branch vertices.

\medskip

Let $s,w$ and $t$ be the three branch vertices in $T$ such that
$d_T(s)=d_T(w)=d_T(t)=3$ and $w\in V(P_T[s,t])$. For each $1\leq
i\leq 5$, let $B_i$ be the vertex set of the connected component of
$T-\{s,w,t\}$ containing $u_i$ and let $v_i$ be the unique vertex in
$B_i\cap N_T(\{s,w,t\})$. Without loss of generality, we may assume
that $\{v_1,v_2\}\subseteq N_T(s)$, $\{v_3,v_4\}\subseteq N_T(t)$
and $v_5\in N_T(w)$. For each $1\leq i\leq 5$ and $x\in B_i$, we use
$x^-$ to denote the predecessor of $x$ on $P_T[s,u_i]$ or
$P_T[t,u_i]$ or $P_T[w,u_i]$. Define $P:=V(P_T[s,t])-\{s,w,t\}$.

For this case, we further choose $T$ such that
\begin{itemize}
\item [$($C3$)$]
$d_T[s,t]$ is as small as possible.
\end{itemize}

It is easy to check that the following claim still holds in this
case. (The proof is exactly the same as that of
Claim~\ref{claim2.2}.)

\begin{claim}\label{claim2.14}
For all $1\leq i,j\leq 5$ and $i\neq j$, if $x\in N(u_j)\cap B_i$,
then $x\notin \{u_i,v_i\}$ and $x^-\notin N(U-\{u_j\})$.
\end{claim}

By applying Claim~\ref{claim2.14}, we deduce that $U$ is an
independent set in $G$.

\begin{claim}\label{claim2.15}
$N(u_i)\cap P=\emptyset$ for each $1\leq i\leq 4$.
\end{claim}

\pf Suppose to the contrary that there exists some vertex $x\in P$
such that $xu_i\in E(G)$ for some $1\leq i\leq 4$. Without loss of
generality, we may assume that $x\in V(P_T[s,w])-\{s,w\}$. Let
$T':=T-v_iv_i^-+xu_i$. If $i\in\{1,2\}$, then $T'$ is a tree in $G$
with $5$ leaves such that $V(T')=V(T)$, $T'$ has three branch
vertices $x,w$ and $t$, $d_{T'}(x)=d_{T'}(w)=d_{T'}(t)=3$, $w\in
V(P_{T'}[x,t])$ and $d_{T'}[x,t]<d_T[s,t]$, contradicting the
condition (C3). Hence we have $i\in\{3,4\}$. Now, $T'$ is a tree in
$G$ with $5$ leaves such that $V(T')=V(T)$, $T'$ has three branch
vertices $s,x$ and $w$, $d_{T'}(s)=d_{T'}(x)=d_{T'}(w)=3$, $x\in
V(P_{T'}[s,w])$ and $d_{T'}[s,w]<d_T[s,t]$. But this also
contradicts the condition (C3). So the claim holds.     \qed

\begin{claim}\label{claim2.16}
$N(u_i)\cap \{w,t\}=\emptyset$ for each $1\leq i\leq 2$ and
$N(u_j)\cap \{s,w\}=\emptyset$ for each $3\leq j\leq 4$.
\end{claim}

\pf Suppose there exists some $i\in \{1,2\}$ such that $wu_i\in
E(G)$ or $tu_i\in E(G)$. Then
$$T':=\left\{\begin{array}{ll}T-sv_i+wu_i,
& \;\mbox{ if } wu_i\in
E(G),\\
T-sv_i+tu_i, & \;\mbox{ if } tu_i\in E(G),
\end{array}\right.$$
is a tree in $G$ with $5$ leaves such that $V(T')=V(T)$ and $T'$ has
two branch vertices $w$ and $t$. By the same argument as in the
proof of Case 1, we can obtain a contradiction. Therefore,
$N(u_i)\cap \{w,t\}=\emptyset$ for each $1\leq i\leq 2$.

By a similar argument as above (by exchanging the roles of $s$ and
$t$), we can also show that $N(u_j)\cap \{s,w\}=\emptyset$ for each
$3\leq j\leq 4$. This completes the proof of Claim~\ref{claim2.16}.
\qed

\begin{claim}\label{claim2.17}
$N(u_5)\cap \{s,t\}=\emptyset$.
\end{claim}

\pf Suppose for a contradiction that $N(u_5)\cap \{s,t\}\neq
\emptyset$. Then we have $su_5\in E(G)$ or $tu_5\in E(G)$. Define
$$T':=\left\{\begin{array}{ll}T-wv_5+su_5,
& \;\mbox{ if } su_5\in
E(G),\\
T-wv_5+tu_5, & \;\mbox{ if } tu_5\in E(G).
\end{array}\right.$$
Now, $T'$ is a tree in $G$ with $5$ leaves such that $V(T')=V(T)$
and $T'$ has two branch vertices $s$ and $t$. By the same argument
as in the proof of Case 1, we can derive a contradiction. So the
assertion of the claim holds.      \qed

\begin{claim}\label{claim2.18}
$N_{\geq 2}(U-\{u_i\})\cap B_i=\emptyset$ for each $1\leq i\leq 5$.
\end{claim}

\pf Suppose the assertion of the claim is false. Then there exists
some vertex $x\in N_{\geq 2}(U-\{u_i\})\cap B_i$ for some $1\leq
i\leq 5$. By Claim~\ref{claim2.14}, we see that
$x\notin\{u_i,v_i\}$. Since $x\in N_{\geq 2}(U-\{u_i\})$, there must
exist two distinct $j,k\in\{1,2,3,4,5\}-\{i\}$ such that
$xu_j,xu_k\in E(G)$. By symmetry between $j$ and $k$, we can always
choose $j$ such that $v_i^-\neq v_j^-$. Let
$T':=T-\{v_iv_i^-,v_jv_j^-\}+\{xu_j,xu_k\}$ and let
$y:=\{s,w,t\}-\{v_i^-,v_j^-\}$. Then $T'$ is a tree in $G$ with $5$
leaves such that $V(T')=V(T)$, $T'$ has two branch vertices $x$ and
$y$, $d_{T'}(x)=4$ and $d_{T'}(y)=3$. By the same argument as in the
proof of Case 1, we can deduce a contradiction. This proves
Claim~\ref{claim2.18}.    \qed

\bigskip

By applying Claim~\ref{claim2.14}, we conclude that $\{u_i\}$,
$N(u_i)\cap B_i$ and $(N(U-\{u_i\}))^-\cap B_i$ are pairwise
disjoint subsets in $B_i$ for each $1\leq i\leq 5$, where
$(N(U-\{u_i\}))^-\cap B_i=\{x^-\;|\;x\in N(U-\{u_i\})\cap B_i\}$. It
follows from Claims~\ref{claim2.15}--\ref{claim2.18} that
$N_5(U)=N_4(U)=N_3(U)=(N_2(U)-N(u_i))\cap B_i=\emptyset$ (for each
$1\leq i\leq 5$). Therefore, for each $1\leq i\leq 5$, we have
\begin{align*}
|B_i|&\geq 1+|N(u_i)\cap B_i|+|(N(U-\{u_i\}))^-\cap B_i|\\
&=
1+|N(u_i)\cap B_i|+|N(U-\{u_i\})\cap B_i|\\
&=1+\sum\limits_{j=1}^5|N(u_j)\cap B_i|.  \tag{$4$}
\end{align*}

By Claims~\ref{claim2.16} and~\ref{claim2.17}, we know that
\begin{align*}
\sum\limits_{i=1}^5|N(u_i)\cap\{s,w,t\}|&=\sum\limits_{i=1}^2|N(u_i)\cap\{s,w,t\}|
+\sum\limits_{i=3}^4|N(u_i)\cap\{s,w,t\}|+|N(u_5)\cap\{s,w,t\}|\\
&\leq 2+2+1\\
&=5.
\end{align*}
On the other hand, by Claim~\ref{claim2.15}, we have
\begin{align*}
\sum\limits_{i=1}^5|N(u_i)\cap P|=|N(u_5)\cap P|\leq |P|.
\end{align*}
Hence
\begin{align*}
|V(P_T[s,t])|=3+|P|=5+|P|-2\geq
\sum\limits_{i=1}^5|N(u_i)\cap\{s,w,t\}|+\sum\limits_{i=1}^5|N(u_i)\cap
P|-2.   \tag{$5$}
\end{align*}

Since $N(U)\subseteq V(T)$ and by (4) and (5), we deduce that
\begin{align*}
|V(T)|&=\sum\limits_{i=1}^{5}|B_i|+|V(P_{T}[s,t])|\\
&\geq \sum\limits_{i=1}^{5}\left(1+\sum\limits_{j=1}^5|N(u_j)\cap
B_i|\right)+\left(\sum\limits_{i=1}^5|N(u_i)\cap\{s,w,t\}|+\sum\limits_{i=1}^5|N(u_i)\cap
P|-2\right)\\
&=5+\sum\limits_{i=1}^{5}\sum\limits_{j=1}^{5}|N(u_j)\cap
B_i|+\sum\limits_{i=1}^{5}|N(u_i)\cap \{s,w,t\}|+\sum\limits_{i=1}^{5}|N(u_i)\cap P|-2\\
&=\sum\limits_{j=1}^{5}|N(u_j)\cap V(T)|+3\\
&=\sum\limits_{j=1}^{5}d(u_j)+3\\
&=d(U)+3.
\end{align*}
This implies that
\begin{align*}
\sigma_5(G)\leq d(U)\leq |V(T)|-3\leq n-3,
\end{align*}
contradicting the assumption that $\sigma_5(G)\geq n-1$. This
completes the proof of Theorem~\ref{theo1.9}. \qed

\bigskip

{\bf Acknowledgements.} The first author was supported by the
National Natural Science Foundation of China (No. 11526160) and the
Science and Technology Innovation Project of Wuhan Textile
University. The second author was supported by the NAFOSTED Grant of
Vietnam (No. 101.04-2018.03).


\begin{thebibliography}{99}
\addtolength{\baselineskip}{-1ex}

\bibitem{BT98}
H. Broersma and H. Tuinstra, Independence trees and Hamilton cycles,
{\it J. Graph Theory} {\bf 29} (1998) 227--237.

\bibitem{CLX17}
X. Chen, M. Li and M. Xu, Spanning $3$-ended trees in $k$-connected
claw-free graphs, {\it Ars Combin.} {\bf 131} (2017) 161--168.

\bibitem{CCH14}
Y. Chen, G. Chen and Z. Hu, Spanning $3$-ended trees in
$k$-connected $K_{1,4}$-free graphs, {\it Sci. China Math.} {\bf 57}
(2014) 1579--1586.

\bibitem{Di05}
R. Diestel, Graph Theory, 3rd Edition, Springer, Berlin, 2005.

\bibitem{GHHSV04}
L. Gargano, M. Hammar, P. Hell, L. Stacho and U. Vaccaro, Spanning
spiders and light-splitting switches, {\it Discrete Math.} {\bf 285}
(2004) 83--95.

\bibitem{KKMOSY12}
M. Kano, A. Kyaw, H. Matsuda, K. Ozeki, A. Saito and T. Yamashita,
Spanning trees with a bounded number of leaves in a claw-free graph,
{\it Ars Combin.} {\bf 103} (2012) 137--154.

\bibitem{Ky09}
A. Kyaw, Spanning trees with at most $3$ leaves in $K_{1,4}$-free
graphs, {\it Discrete Math.} {\bf 309} (2009) 6146--6148.

\bibitem{Ky11}
A. Kyaw, Spanning trees with at most $k$ leaves in $K_{1,4}$-free
graphs, {\it Discrete Math.} {\bf 311} (2011) 2135--2142.

\bibitem{LV71}
M. Las Vergnas, Sur une propriet\'{e} des arbres maximaux dans un
graphe, {\it C. R. Acad. Sci. Paris Ser. A} {\bf 272} (1971)
1297--1300.

\bibitem{MOY14}
H. Matsuda, K. Ozeki and T. Yamashita, Spanning trees with a bounded
number of branch vertices in a claw-free graph, {\it Graphs Combin.}
{\bf 30} (2014) 429--437.

\bibitem{MS84}
M. M. Matthews and D. P. Sumner, Hamiltonian results in
$K_{1,3}$-free graphs, {\it J. Graph Theory} {\bf 8} (1984)
139--146.

\bibitem{OY11}
K. Ozeki and T. Yamashita, Spanning trees: A survey, {\it Graphs
Combin.} {\bf 27} (2011) 1--26.

\bibitem{Wi79}
S. Win, On a conjecture of Las Vergnas concerning certain spanning
trees in graphs, {\it Results Math.} {\bf 2} (1979) 215--224.

\end{thebibliography}
\end{document}